\documentclass[11pt,draft]{article}
\usepackage{amsfonts,amssymb,amsmath,amsthm}
\usepackage{parskip}
\newtheorem{definition}{Definition}%[section]
\newtheorem{proposition}[definition]{Proposition}%[section]

\numberwithin{equation}{section}
\numberwithin{definition}{section}
\newcommand{\Ud}{\mathrm{d}}
\newcommand{\Iff}{if\textcompwordmark f}
\newcommand{\rR}{\mathbb{R}}
\DeclareMathOperator{\rank}{rank}
\DeclareMathOperator{\divz}{div}
\DeclareMathOperator{\sign}{sign}

\begin{document}

\title{Notes on relative normalizations of ruled surfaces in the three-dimensional Euclidean space}

\author{{Stylianos Stamatakis and Ioanna-Iris Papadopoulou}\textbf{\medskip}\\ \emph{Aristotle University of Thessaloniki}\\ \emph{Department of Mathematics}\\ \emph{GR-54124 Thessaloniki, Greece}\\  \emph{e-mail: stamata@math.auth.gr}}
\date{}
\maketitle
\begin{abstract}
\noindent This paper deals with relative normalizations of skew ruled surfaces in the Euclidean space $\mathbb{E}^{3}$. In section 2 we investigate some new formulae concerning the Pick invariant, the relative curvature, the relative mean curvature and the curvature of the relative metric of a relatively normalized ruled surface $\varPhi$ and in section 3 we introduce some special normalizations of it. All ruled surfaces and their corresponding normalizations that make $\varPhi$ an improper or a proper relative sphere are determined in section 4. In the last section we study ruled surfaces, which are \emph{centrally} normalized, i.e., their relative normals at each point lie on the corresponding central plane. Especially we study various properties of the Tchebychev vector field.  We conclude the paper by the study of the central image of $\varPhi$.


\end{abstract}
%\vskip 1 cm

%\begin{center}
%Dedicated to on the occasion of his  $^{th}$ birthday
%\end{center}
\section{Preliminaries}\label{Sec1}
%This paper deals with skew ruled surfaces $\varPhi$ in the Euclidean space $\mathbb{E}^{3}$ which are relative normalized and is organized as follows: After the introductory Sect. 1, we investigate in Sect. 2 some new formulae concerning the Pick invariant, the relative curvature, the relative mean curvature and the curvature of the relative metric of a relatively normalized ruled surface.
%In Sect. 3 we introduce some special normalizations for a given ruled surface. In Sect. 4 we determine all ruled surfaces and their corresponding normalizations that make $\varPhi$ an improper or a proper relative sphere. In Sect. 5 we study ruled surfaces centrally normalized. In Sect. 6 we find the relative image of a centrally normalized ruled surface, in Sect. 7 we study some properties of the Tchebychev vector field and in Sect. 8 we we study ruled surfaces polarly normalized.

To set the stage for our work we present briefly some elementary facts regarding the relative Differential Geometry of surfaces and the Differential Geometry of ruled surfaces in the Euclidean space $\mathbb{E}^{3}$; we shall follow the notations and definitions of \cite{Pottmann} and \cite{Schirokow}.

In the three-dimensional Euclidean space $\mathbb{E}^{3}$ let $\varPhi$ be a ruled $C^{r}$-surface of nonvanishing Gaussian curvature, $r\geq3$, defined by an injective $C^{r} $-immersion $\overline{x}=\overline{x}(u,v)$  on a region $U:=I\times \rR \,(I\subset \rR$  open interval)  of $\rR^{2}$.
Let $\langle \,,\rangle$ denote the standard scalar product in $\mathbb{E}^{3}$.
We introduce the so-called standard parameters $u\in I, v \in \rR$ of $\varPhi$, such that
    \begin{equation}    \label{1}
        \overline{x}(u,v)=\overline{s}(u) + v\,\overline{e}(u),
    \end{equation}%
with%
    \begin{equation}    \label{5}
        \left \vert \overline{e}΄\right \vert =|\overline{e}'|=1,\quad \langle \overline{s}',\overline{e}'\rangle =0,
    \end{equation}%
where the differentiation with respect to $u $ is denoted by a prime.
Here $\varGamma: \overline{s}=\overline{s}(u)$ is the striction curve of $\varPhi$ and the parameter $u$ is the arc length along the spherical curve $\overline{e}=\overline{e}(u)$.

Let $$\delta(u) :=(\overline{s}',\overline{e},\overline{e}')$$ be the distribution parameter, $$\kappa(u):=(\overline{e},\overline{e}',
\overline{e}'')$$ be the conical curvature and $$\sigma(u) :=\sphericalangle (\overline{e},\overline{s}'), \quad \text{where} \quad -\frac{\pi }{2}<\sigma \leq \frac{\pi }{2}, \,\,\, \sign \sigma = \sign \delta)$$ be the striction of $\varPhi$. %are the fundamental invariants of $\varPhi$ and determine uniquely, up to Euclidean rigid motions, the ruled surface $\varPhi $.
We consider yet the central normal vector $\overline{n}(u):=\overline{e}'$ and the central tangent vector $\overline{z}(u):=\overline{e} \times \overline{n}$. %The moving frame $\mathcal{D}  $ of $\varPhi $ is the orthonormal frame which is attached to
The moving frame $\mathcal{D} : = \{\overline{e}, \overline{n}, \overline{z}\}$ of $\varPhi$ %the striction point $\overline{s}(u)$, and consists of the vector $\overline{e}(u) $, It
fulfils the equations \cite[p. 280]{Pottmann}
    \begin{equation}\label{10}
    \overline{e}'=\overline{n},\quad \overline{n}'=-\overline{e}+\kappa \,\overline{z},\quad \overline{z}'=-\kappa \,\overline{n}.
    \end{equation}
Then, we have
    \begin{equation}    \label{15}
    \overline{s}'=\delta \,\lambda \, \overline{e}+\delta \, \overline{z}, \quad \text{where}\quad  \lambda(u) :=\cot \sigma.
    \end{equation}%
We denote partial derivatives of a function (or a vector-valued function) $f$ in the coordinates $u^{1}:=u,\,u^{2}:=v$ by
    %$$
    %f_{/i}:=\frac{\partial f}{\partial u^{i}},~f_{/ij}:=\frac{\partial ^{2}f}{\partial u^{i}\partial u^{j}} \quad \text{etc.}
    %$$
$f_{/i}, f_{/ij}$ etc.
Then from (\ref{1}) and (\ref{15})\ we obtain
    \begin{equation}\label{20}
    \overline{x}_{/1}=\delta \, \lambda \,\overline{e}+v\,\overline{n}+\delta \,\overline{z},\quad \overline{x}_{/2}=\overline{e},
    \end{equation}
and thus the unit normal vector $\overline{\xi}(u,v) $ to $\varPhi $ is expressed by
    \begin{equation}\label{25}
    \overline{\xi}=\frac{\delta \,\overline{n}-v\,\overline{z}}{w},\quad \text{where}\quad w:=\sqrt{\delta ^{2} + v^{2}}.
    \end{equation}
Let $I=g_{ij} \Ud u^i \Ud u^j$ be the first and $II=h_{ij}\Ud u^i \Ud u^j$ be the second fundamental form of $\varPhi $, where
    \begin{equation}\label{30}
         g_{11}=w^{2}+\delta ^{2}\,\lambda ^{2},\qquad g_{12}=\delta\, \lambda, \qquad   g_{22}=1,
    \end{equation}
    \begin{equation}\label{35}
     h_{11}=-\frac{\kappa\, w^{2}+\delta ' \,v-\delta ^{2}\,\lambda }{w},\qquad h_{12}=\frac{\delta }{w}, \qquad  h_{22}=0.
    \end{equation}
The Gaussian curvature $\widetilde{K}(u,v)$ and the mean curvature $H_I(u,v)$ of $\varPhi $ are respectively given by \cite{Pottmann}
   \begin{equation}\label{40}
  \widetilde{K}=- \frac{\delta ^{2}}{w^{4}}, \quad H_{I}=-\frac{\kappa w^{2}+\delta ' v+\delta^{2}\lambda}{2w^{3}}.
 \end{equation}
%and $H_{I}$ by
 %  \begin{equation}\label{41}
% \end{equation}

A $C^{s}$-relative normalization of $\varPhi$ is a $C^{s}$-mapping $\overline{y} = \overline{y}(u,v), 1\leq s < r$, defined on $U$, such that
    \begin{equation}\label{45}
    \rank (\{\overline{x}_{/1},\overline{x}_{/2},\overline{y}\})=3,\,\,
    \rank (\{\overline{x}_{/1},\overline{x}_{/2},\overline{y}_{/i}\})=2,\,\,
    i=1,2,\,\,\forall \left(u,v\right) \in U.
    \end{equation}
The pair $\left(\varPhi,\overline {y}\right)$ is called a relatively normalized ruled surface and the line issuing from a point $P \in \varPhi$ in the direction $\overline{y}$ is called the relative normal of $\varPhi$ at $P$. When we move the vectors $\overline{y} $
%of the relative normalization
to the origin, the endpoints of them describe the relative image of $\varPhi$.

Let $q(u,v):=\langle \bar{\xi},\bar{y}\rangle$, denote the support function of the relative normalization $\bar{y} $ (see~ \cite{Manhart3}). As follows from (\ref{45}) $q $ never vanishes on $U $.
Conversely, when a support function $q$ is given, the relative normalization $\overline{y}$ of the ruled surface $\varPhi$ is uniquely determined and can be expressed in terms of the moving frame $\mathcal{D}$ as follows \cite[p.179]{Stamatakis3}:
    \begin{equation}\label{70}
    \overline{y}=y_{1}\,\overline{e}+y_{2}\,\overline{n}+y_{3}\,\overline{z},
    \end{equation}
where
    \begin{equation}\label{75}
    y_{1}=-w\frac{\delta q_{/1}+q_{/2}(\kappa \,w^{2}+\delta '\,v)}{\delta ^{2}},\quad
    y_{2}=\frac{\delta ^{2}\,q-w^{2}\,v\,q_{/2}}{\delta w},\quad
    y_{3}=-\frac{v\,q+w^{2}\,q_{/2}}{w}.
    \end{equation}
One can easily verify the following relations:
    \begin{equation}\label{76}
    y_{1}+y_{2/1}-\kappa y_{3}=\frac{v}{\delta}(y_{3/1}+\kappa y_{2}),\quad \quad y_{2/2}=\frac{v}{\delta}y_{3/2}.
    \end{equation}
For the coefficients $G_{ij}$ of the relative metric $G$ of $(\varPhi, \overline{y})$, which is indefined, the following applies
    \begin{equation}\label{51}
    G_{ij}=q^{-1}\,h_{ij}.
    \end{equation}
%Furthermore, on account of (\ref{35}), the coefficients of the  inverse relative metric tensor are computed by
%The covector $\bar{X} $ of the tangent plane is defined by
   % \begin{equation}\label{50}
    %\langle \bar{X},\bar{x}_{/i}\rangle =0 \quad (i=1,2)\quad \text{and}\quad \langle \bar{X},\bar{y}\rangle =1
    %\end{equation}%
%and the  relative metric $G=G_{ij}\Ud u^i \Ud u^j$ by $G_{ij}:=\langle \bar{X},\bar{x}_{/ij}\rangle.$
%One can easily cf. that $\bar{X}=q^{-1}\bar{\xi} $ and
%\begin{equation}\label{51}
%G_{ij}=q^{-1}\,h_{ij},
%\end{equation}
%where $q:=\langle \bar{\xi},\bar{y}\rangle$ denotes the \emph{support function of the relative normalization} $\bar{y} $ (cf.
Then, on account of (\ref{35}), the coefficients of the  inverse relative metric tensor are computed by
    \begin{equation}\label{52}
    G^{(11)}=0,\quad G^{(12)}=\frac{w\,q}{\delta },\quad G^{(22)}= w\,q \, \frac{\kappa\, w^{2}+\delta '\,v-\delta ^{2}\,\lambda }{\delta ^{2}}.
    \end{equation}
%It is well known that, when a support function $q(u,v) $ is given, then the relative normalization $\overline{y} $ is uniquely determined by (cf.~
%\cite{Manhart3})
 %   \begin{equation}\label{55}
  %  \overline{y}=-h^{\left( ij\right) }\,\,q_{/i}\,\,\overline{x}_{/j}+q\,\overline{\xi},
   % \end{equation}%
%where $h^{(ij)} $ are the coefficients of the inverse tensor of $h_{ij} $.
%By using this last equation, (\ref{20}), (\ref{25}) and (\ref{35}), the relative normalization $\overline{y}$ of $\varPhi $ can be expressed in terms of the vectors of the moving frame $\mathcal{D} $, as follows:
 %   \begin{equation}\label{70}
  %  \overline{y}=y_{1}\,\overline{e}+y_{2}\,\overline{n}+y_{3}\,\overline{z},
   % \end{equation}
%where
 %   \begin{equation}\label{75}
  %  y_{1}=-w\frac{\delta q_{/1}+q_{/2}(\kappa \,w^{2}+\delta '\,v)}{\delta ^{2}},\quad
   % y_{2}=\frac{\delta ^{2}\,q-w^{2}\,v\,q_{/2}}{\delta w},\quad
    %y_{3}=-\frac{v\,q+w^{2}\,q_{/2}}{w}.
    %\end{equation}
For a function (or a vector-valued function) $f$ we denote by $\nabla ^{G}f$ the first Beltrami differential operator and by $\nabla _{i}^{G}f $ the
covariant derivative, both with respect to the relative metric.
The coefficients $A_{ijk}(u,v)$ of the Darboux tensor are given by
    \begin{equation}\label{90}
    A_{ijk}:= \frac{1}{q} \, \langle \overline{\xi},\,\nabla _{k}^{G}\,\nabla _{j}^{G}\,\overline{x}_{/i}\rangle.
    \end{equation}
Then, by using the relative metric tensor $G_{ij} $ for ``raising and lowering
the indices'', the Pick invariant $J(u,v)$ of $(\varPhi,\overline{y})$ is defined by
    \begin{equation*}
      J:=\frac{1}{2}A_{ijk}\,A^{ijk}.
    \end{equation*}
%We denote by $q_{AFF}=|\widetilde{K}|^{1/4}$ the support function of the equiaffine normalization $\overline{y}_{AFF}$.
In \cite{Stamatakis3} it was shown, that the coefficients of the Tchebychev vector $\overline{T}(u,v)$ of $(\varPhi,\overline{y})$, which is defined by
    \begin{equation*}
    \overline{T}:=T^{m}\, \overline{x}_{/m},\quad \text{where\quad }T^{m}:=\frac{1}{2}A_{i}^{im},
    \end{equation*}
are given by
   % \begin{equation*}
    %\overline{T}=\bigtriangledown^{G}(\ln\phi, \overline{x}),
    %\end{equation*}
%where $\varphi=\frac{q}{q_{AFF}}$ is the Tchebychev function of $\overline{y}$ and in \cite{Stamatakis3} it was shown that in the case of a relatively normalized ruled surface $\varPhi$ holds
    \begin{equation}\label{121}
    T^{1}=\frac{w^{2}q_{/2}+vq}{\delta w},\,\, T^{2}=\frac{2\delta\,w^{2}q_{/1}+\delta'\,q(\delta^{2}-v^{2})}{2\delta^{2}\, w}+\frac{T^{1}(\kappa w^{2}+\delta'\,v-\delta^{2}\lambda)}{\delta}.
    \end{equation}
$\overline{T}$ can be expressed in terms of the moving frame $\mathcal{D}$ as follows \cite{Stamatakis3}
    \begin{equation}\label{120}
    \overline{T}=w\frac{q\left( 2\kappa v+\delta ^{\prime }\right) +2\delta
    q_{/1}+2q_{/2}(\kappa w^{2}+\delta ^{\prime }v)}{2\delta ^{2}}\bar{e}+\frac{%
    vq+w^{2}q_{/2}}{\delta w}\left( v\bar{n}+\delta \bar{z}\right).
    \end{equation}
%To compute its components it is convenient to use the well known relations
%    \begin{equation}\label{120}
 %   T^{i}=\left[ \ln \left( \frac{/q|}{q_{AFF}}\right) \right] _{/j}G^{(ij)},
  %  \end{equation}%
%where, by virtue of (\ref{40}),
 %   \begin{equation}\label{125}
  %  q_{AFF}=|\widetilde{K}|^{1/4}=|\delta|^{1/2}\,w^{-1}
   % \end{equation}
%denotes the support function of the equiaffine normalization $\overline{y}_{AFF}$.
The relative shape operator has the coefficients $B_{i}^{j}(u,v)$ defined by
    \begin{equation}\label{125}
    \overline{y}_{/i}=:-B_{i}^{j}\, \overline{x}_{/j}.
    \end{equation}
Then, the relative curvature $K(u,v)$ and the relative mean curvature $H(u,v)$ are defined by
    \begin{equation}\label{130}
    K:=\det \left(B_{i}^{j}\right),\quad H:=\frac{B_{1}^{1}+B_{2}^{2}}{2}.
    \end{equation}
%Let $S$ be the scalar curvature of the relative metric $G$, which is defined formally as the
%curvature of the pseudo-Riemannian manifold ($\varPhi,G$). It is known, that the following identity holds (\emph{Theorema Egregium of the relative Differential Geometry}, cf. \cite[p. 197]{Manhart3})
    %\begin{equation*}
    %J+H-S=2\|\overline{T} \|_{G},
    %\end{equation*}
%where
   % \begin{equation}\label{134}
    %\|\overline{T} \|_{G}=G_{ij}\,T^{i}\,T^{j}
    %\end{equation}
%$\|\overline{T} \|_{G}$ is the norm of the Tchebychev vector with respect to the relative metric.
We mention finally for later use, that among the surfaces $\varPhi \subset \mathbb{E}^3$ with negative Gaussian curvature the ruled surfaces are characterized by the relation \cite{Stamatakis4}
    \begin{equation}\label{136}
    3 H - J -3 S = 0,
    \end{equation}
where $S(u,v)$ is the scalar curvature of the relative metric $G$, which is defined formally as the curvature of the pseudo-Riemannian manifold ($\varPhi,G$).
%By using (\ref{35}), (\ref{50}), (\ref{65b}), (\ref{110}) and (\ref{130}) we find
%\begin{equation}\label{144}
 % \|\overline{T}\|_{G}=\frac{2}{3}J.
%\end{equation}
%So, by taking (\ref{65a}) into consideration we arrive to the following
%\begin{proposition}
%For each ruled surface $\varPhi \subset \mathbb{R}$ holds
%We mention finally for later use the relation
%\begin{equation}\label{144a}
%y_{3/2}\left( y_{1/1}-y_2 \right) - y_{1/2}\left(  \kappa \, y_2+y_{3/1} \right) = -\delta \, K.
%\end{equation}

\section{Some formulae for $J, K, H$ and $S$}\label{Sec2}
In this section we express the relative magnitudes $J,K,H$ and $S$ of the relatively normalized ruled surface $\varPhi$ in terms of the fundamental invariants $\delta,\kappa$ and $ \lambda$ of $\varPhi$ and the support function $q$. Firstly we compute the Pick invariant $J$. We notice that by virtue of the symmetry of the Darboux tensor \eqref{90} we have
    \begin{equation}\label{200}
    J=\frac{3}{2} \left(A_{112} A^{112}+A_{122} A^{122} \right)+\frac{1}{2}\left(A_{111} A^{111}+A_{222} A^{222}\right).
    \end{equation}
By using the well known equation \cite[p. 196]{Manhart3}
    \begin{equation*}
    A_{ijk}=\frac{1}{q}\langle \overline{\xi},\overline{x}_{/ijk}\rangle -\frac{1}{2}%
    \left( G_{ij/k}+G_{jk/i}+G_{ki/j}\right) ,
    \end{equation*}
and the relations  \eqref{10}, \eqref{20}, \eqref{25}, \eqref{35}, \eqref{51} and \eqref{52} we get by straightforward calculations
%\begin{subequations}%\label{A}
    \begin{align*}
    A^{111}&=A_{222}=0, \\
    A_{112}&=\frac{-1}{2q^2\, w^3} \big\{ \left(w^2 \,q_{/2}+q\,v \right) \big[\kappa\, v^2+ \delta '\,v+\delta^2 \, (\kappa-\lambda)\big]+\delta'\,w^2 q-2\delta\,(w^2 \,q_{/1} + \delta \,\delta '\, q)\big\}, \\
    A^{112} & = \frac{q}{\delta^2}  \left(w^2\, q_{/2}+q\,v\right),\\
    A_{122}&= \frac{\delta}{w^3 \,q^2} \left(w^2 \,q_{/2}+q\,v\right),\\
    A^{122}&=\frac{q}{2 \delta ^3}\big\{3 \left(w^2\, q_{/2}+q\,v \right) \big[\kappa \,v^2+ \delta '\,v+\delta^2\, (\kappa-\lambda)\big]-\delta'\,w^2 \,q+2\delta(w^2\, q_{/1} + \delta \,\delta '\, q)\big\}.
    \end{align*}
%\end{subequations}
Inserting these relations in \eqref{200} we obtain
   \begin{equation}\label{210}
    J = \frac{3\left(w^2 q_{/2}+v\,q\right)}{2\delta^2 w^3 \,q}\Big\{w^2 \!\left[ \kappa  q v + 2 \delta q_{/1} +q_{/2} \left(\kappa \,w^2 + \delta ' \, v-\delta^2 \lambda \right)  \right]  -\delta^2 q \left( \lambda v - \delta '  \right) \Big\}.
    \end{equation}
Next, we wish now to compute the relative curvature and the relative mean curvature. To this end we find  on account of \eqref{10}, \eqref{20},  \eqref{70},  \eqref{75} and \eqref{125} firstly the coefficients $B_i^j$ of the relative shape operator:
%\begin{subequations}%\label{B}
    \begin{align*}
    B_{1}^{1}&= -\frac{1}{\delta^{2} w^{3}}\big[\delta^{2}\delta' q v + \kappa w^{2} \left(\delta^{2}q - v w^{2} q_{/2}\right)-\delta w^{2} \left(\delta\delta' q_{/2} + v q_{/1} + w^{2} q_{/12}\right)\big],\\
        \begin{split}
        B_{1}^{2} & =\frac{1}{\delta^{3}w^{3}}\big\{\delta^{4}w^{2}q - \delta^{2}v w^{4}q_{/2} -\delta' w^{2}\left(\delta^{2}+2v^{2}\right)\big[q_{/2}\left(\kappa w^{2}+\delta' v\right)+\delta q_{/1}\big]\\
        & + \delta^{2} \lambda \big[\delta^{2}\delta'\left(q v + w^{2}q_{/2}\right) + \kappa w^{2}\left(\delta^{2}q - v w^{2} q_{/2}\right) -\delta  w^{2}\left(2 \delta \delta' q_{/2} + q_{/1} v + w^{2}q_{/12}\right)\big] \\
        & +\delta w^{4}\big[q_{/2}\left(2 \delta \delta' \kappa + \kappa' w^{2} + \delta'' v\right) + \delta' q_{/1} + q_{/12} \left(\kappa w^{2} + \delta' v\right) + \delta q_{/11}\big]\big\},
        \end{split}\\
     B_{2}^{1}&= \frac{1}{\delta w^{3}}\left(2 q_{/2}v w^{2} + \delta^{2}q + w^{4} q_{/22}\right),\\
        \begin{split}
        B_{2}^{2}&= \frac{1}{\delta^{2}w^{3}}\big\{-\delta^{2} \lambda\left(2 q_{/2}v w^{2}+ \delta^{2}q + w^{4}q_{/22}\right)+v w^{2}\big[q_{/2} \left(\kappa w^{2} + \delta' v\right)+ \delta q_{/1}\big] \\
        & + w^{4}\big[q_{/2} \left(2 \kappa v + \delta'\right)+q_{/22} \left(\kappa w^{2}+\delta' v\right)+ \delta q_{/12}\big]\big\}.
        \end{split}
    \end{align*}
%\end{subequations}
Substituting the above relations in \eqref{130} we get
     \begin{equation}\label{215}
    \begin{split}
    K  =\frac{-1}{\delta^4 \, w^6} & \Big\{ \big\{ \lambda \, \delta^2 [ -2v^3 q_{/2} -\delta^2 ( q+2 v q_{/2})-w^4 \, q_{/22} ]  \\
    & +v\,w^2[ q_{/2} ( \kappa w^2 + \delta' v ) + \delta \, q_{/1} ] +w^4 [q_{/2} ( 2\kappa \, v+\delta' )+q_{/22} ( \kappa  w^2 +\delta' \, v )    \\
    &+\delta q_{/12}] \big\}[ \delta^2 \, \delta' ( q\, v + w^2 \, q_{/2} )+ \kappa w^2 ( \delta^2 \, q - v \, w^2 \, q_{/2} )   \\
    & -\delta \, w^2 ( 2\delta\, \delta' \, q_{/2}+ v\, q_{/1} +w^2 q_{/12})] + [ -2q_{/2} \, v^3 -\delta^2(q+2 v q_{/2})  \\
    & -w^4 \, q_{/22}] \big\{\delta^4 \, w^2 \,q-\delta^2 \, v \,w^4 \,q_{/2}+ \delta^2 \,\delta' w^2 [ q_{/2} ( \kappa w^2+\delta' v  )  +\delta q_{/1} ]    \\
    &-2 \delta' w^4 [q_{/2} (\kappa w^2 + \delta' v )+\delta q_{/1}]+\delta \lambda [ \delta ^3 \delta' (q v +w^2 q_{/2}) \\
    &+ \delta \kappa w^2 ( \delta^2 q -v w^2 q_{/2} )-\delta^2 w ^2( 2\delta \delta' q_{/2} + v q_{/1} + w^2 q_{/12})] \\
    &+ w^4 [ q_{/2} ( 2\delta \delta' \kappa + \kappa' w^2 + \delta'' v )+ \delta' q_{/1} + \delta \, q_{/11} + q_{/12} ( \kappa w^2 + \delta' v) ] \big\}\Big\},
    \end{split}
    \end{equation}
    and
    \begin{equation}\label{220}
    \begin{split}
        H &= \frac{1}{2\delta^{2}w^{3}} \Big \{-\delta^{2}q \left(\kappa w^{2} + \delta' v + \delta^{2} \lambda \right) + 2 w^{2}q_{/2} \big[\left(2\kappa v + \delta'\right)w^{2}-\delta^{2}\lambda v\big]  \\
        &+ w^{4}q_{/22} \left(\kappa w^{2} + \delta' v - \delta^{2} \lambda\right) + 2 \delta v w^{2} q_{/1} + 2 \delta w^{4} q_{/12} \Big \}.
        \end{split}
    \end{equation}
Inserting \eqref{200} and \eqref{220} in \eqref{136} we infer the scalar curvature of the relative metric $G$
    \begin{equation}\label{225}
        \begin{split}
        S& =\frac{1}{2 \delta^{2} w^{3} q}\Big \{q^{2} \big[-\kappa w^{4} + \delta^{2} \big(\lambda v^{2} - 2 \delta' v - \delta^{2}\lambda \big )\big] + w^{4}q\, q _{/2} \big(2 \kappa v + \delta' \big) \\
        & + w^{4} \big(\kappa \, w^{2} + \delta' \, v - \delta^{2} \lambda \big) \big( q\, q_{/22} - q_{/2}^{2}\big)- 2\delta w^{4} q_{/1} q_{/2} + 2 \delta w^{4} q\, q_{/12} \Big \}.
        \end{split}
    \end{equation}
%By means of \eqref{121} we obtain finally the relation
%\begin{equation}\label{240}
 %   \divz ^{G}\overline{T}-2S=\frac{2\,\kappa\,w\, q}{\delta^{2}},
%\end{equation}
%where
%\begin{equation}\label{230}
   % \divz ^{G}\overline{T}=\frac{(|G|^{1/2}\,T^{i})_{/i}}{|G|^{1/2}}, \quad G : = \det \left( G_{ij} \right) = -\frac {\delta^{2}}{w^2 \,q^2},
 %\end{equation}
%is the divergence  of $\overline{T}$ with respect to the relative metric $G$ of $\varPhi$.
The divergence $\divz ^{G}\overline{T}$ of $\overline{T}$ with respect to the relative metric $G$ of $\varPhi$ is given by \cite{Strubecker}
    \begin{equation}\label{230}
    \divz ^{G}\overline{T}=\frac{\left(|G|^{1/2}\,T^{i}\right)_{/i}}{|G|^{1/2}},  \quad \text{where}\quad G : = \det \left( G_{ij} \right) = -\frac {\delta^{2}}{w^2 \,q^2}.
    \end{equation}
%where
%    \begin{equation}\label{235}
 %   G : = \det \left( G_{ij} \right) = -\frac {\delta^{2}}{w^2 \,q^2}.
  %  \end{equation}
By taking \eqref{121}, \eqref{225} and \eqref{230} into consideration it turns out that
    \begin{equation}\label{240}
    \divz ^{G}\overline{T}-2S=\frac{2\,\kappa\,w\, q}{\delta^{2}}.
    \end{equation}

\section{Special relative normalizations}\label{Sec3}
In \cite{Stamatakis3} I.~Kaffas and S.~Stamatakis  have studied the so called asymptotic relative normalizations of a given ruled surface, that is relative normalizations such that the relative normal at each point $P$ of $\varPhi$ lies on the corresponding asymptotic plane $\{P;\overline{e},\overline{n}\}$ of $\varPhi$. Following this idea we consider relative normalizations such that the relative normal at each point $P$ lies a) on the corresponding central plane $\{P;\overline{e},\overline{z}\}$, or b) on the corresponding polar plane $\{P;\overline{n},\overline{z}\}$.

The first case occurs \Iff{}
%The relative normal at each point $P$ of $\varPhi$ lies on the corresponding central plane $\{P;\overline{e}
%,\overline{z}\}$ \Iff{}
 $y_2=0$, or, because of (\ref{75}b), \Iff{} the support function of $\overline{y}$ is of the form
    \begin{equation}\label{301}
    q=\frac{g\,v}{w},
    \end{equation}
where $g=g(u)$ is an arbitrary nonvanishing $C^3$-function. We call the corresponding relative normalization \textit{central}.
Obviously in this case it is
   \begin{equation}\label{305}
    \overline{y}=-\frac{g'\,v+\delta \,\kappa \, g }{\delta}\overline{e}-g\,\overline{z},
    \end{equation}
cf. \eqref{70}, \eqref{75}. The second case occurs \Iff{}
%The relative normal at each point $P$ of $\varPhi$ lies on the corresponding polar plane $\{P;\overline{n}
%,\overline{z}\}$ \Iff{}
 $y_1=0$, or, because of (\ref{75}a), \Iff{} the support function of $\overline{y}$ is of the form
    \begin{equation*}
    q=f(V), \quad \text{where} \quad V =\arctan \frac{v}{\delta}-\int \!\! \kappa \, \Ud u
    \end{equation*}
and $f$ is an arbitrary nonvanishing $C^2$-function of $V$. We call the arising relative normalization \textit{polar}. We find
    \begin{equation*}
    \overline{y}=\frac{\delta \, q-v\,\dot{q}}{w}\overline{n}-\frac{v \, q+\delta \, \dot{q}}{w}\overline{z},
    \end{equation*}
where the dot denotes the derivative in $V$.

Finally, let the relative image be as well as $\varPhi$  a ruled surface whose generators are parallel to those of $\varPhi$. Then  $y_{2/2}=y_{3/2} = 0$, from which, by means of \eqref{75}, we obtain
    \begin{equation*}
    %y_{3/2}=-\frac{2 v \, w^2 \, q_{/2}+\delta^2 \, q + q_{/22} \, w^4}{w^3},
    %\end{equation}
2 v \, w^2 \, q_{/2}+\delta^2 \, q + q_{/22} \, w^4=0.
    \end{equation*}
    Consequently \begin{equation}\label{320}
    q=\frac{f+g\,v}{w},
    \end{equation}
where $f$ and $g$ are arbitrary $C^3$-functions of $u$, such that $q \neq 0$. In this case we have
    \begin{equation}\label{325}
    \overline{y} = \overline{s}^* - \frac {\delta \, g' - \kappa \, f} {\delta ^ 2}v \, \overline{e},
     \end{equation}
where
    \begin{equation*}
    \overline{s}^*(u) = - \bigg[ \left( \frac {f}{\delta} \right) ^{'} + \kappa \, g \bigg] \overline{e} + \frac{f}{\delta} \, \overline{n} - g\,\overline{z}.
    \end{equation*}
From \eqref{325} it follows that the relative image of $\varPhi$ is a curve or a ruled surface whose generators are parallel to those of $\varPhi$ \Iff{} the function $\delta \, g'-\kappa \,f $ vanishes everywhere or nowhere in $I$, respectively. We call in the sequel such a normalization \textit{right}. We recognise immediately that both asymptotic and central normalizations belong to the right ones. In section 5 of this paper we investigate the central normalizations leaving the study of the polar and the right ones for a subsequent paper.

\section{$\varPhi$ is an improper or a proper relative sphere}\label{Sec4}

In this section we investigate all  ruled surfaces $\varPhi$ and the corresponding support functions $q(u,v)$ so that $\varPhi$ is an improper or a proper relative sphere.

It is easily verified from \eqref{10}, \eqref{70} and \eqref{76} that $\varPhi$ is an improper relative sphere, i.e., by definition \cite{Manhart2}, its relative image degenerates into a point
%    \begin{equation} \label{180}
 ($   \overline{y}_{/i}=\overline{0},i=1,2$),
%    \end{equation}
%Firstly, taking on account the relations \eqref{125}, we observe that $ B_{i}^{j}=0$ $\forall \, i,j=1,2$. Consequently, by \eqref{130}, it is $H=K=0$.%
\Iff{} the following relations hold true
    \begin{equation}\label{400}
    y_{1/1} - y_{2} = \kappa\, y_{2} + y_{3/1} = y_{1/2} = y_{3/2} = 0.
    \end{equation}
It is however obvious that $K=H=0$.
By means of $ y_{3/2} = 0$ and (\ref{75}c) we derive that the support function has the form \eqref{320}, i.e., the normalization is right.
We distinguish now two cases:

\textsc{\textbf{Case I}.} $\varPhi$ is conoidal ($\kappa =0 $). From \eqref{75}, \eqref{320} and \eqref{400} we find
    \begin{equation}\label{405}
    f=\delta \left( c_{1} \cos u + c_{2} \sin u \right), \quad g=c_{3}, \quad c_{1},c_{2},c_{3} \in \rR,\quad c_{1}^{2}+c_{2}^{2}+c_{3}^{2}\neq 0.
    \end{equation}
The corresponding relative normalization of $\varPhi$ then results
    \begin{equation}\label{410}
    \overline{y}=\left(c_{1}\,\sin u -c_{2} \cos u \right)\overline{e}+\left(c_{1}\,\cos u +c_{2} \sin u \right)\overline{n}+\overline{c},
    \end{equation}
where $\overline{c} =-c_{3} \, \overline{z}$ is a constant vector. One can easily verify that the converse is valid as well, i.e., that the relative normalization \eqref{410} is constant.
From \eqref{136}, \eqref{210},  \eqref{320} and \eqref{405} we find
    $$
    J = -3S = 3c_{3} \frac { \delta \left( 2c_{2} \cos u - 2c_{1} \sin u - c_{3} \lambda \right) + \delta' \left( c_{1} \cos u + c_{2} \sin u \right) } {2\,\delta \big [ c_{3} v + \delta \left( c_{1} \cos u + c_{2} \sin u   \right)\big ]}.
    $$
\textsc{\textbf{Case II}.} $\varPhi$ is non-conoidal. From the relations \eqref{75} and \eqref{400} we take
    \begin{equation}\label{415}
    f=\frac{\delta \,g'}{\kappa},
    \end{equation}
while the function $g$ fulfils the equation
    \begin{equation}\label{420}
    \left( \frac{g'}{\kappa}\right)^{''} +\frac{g'}{\kappa}+\left( \kappa \,g \right)^{'}=0.
    \end{equation}
In this case we find
    \begin{equation}\label{425}
    \overline{y} = - \bigg [ \kappa \, g + \left( \frac{g'}{\kappa}\right)^{'} \bigg ] \overline{e} + \frac{g'}{\kappa}\, \overline{n} - g \, \overline{z}.
    \end{equation}
The inverse is valid as well: The relative normalization \eqref{425}, under the assumption \eqref{420}, is constant.
From \eqref{136}, \eqref{210},  \eqref{320} and \eqref{415} we obtain
    $$
    J=-3S=3g \, \frac{ \kappa^{2} g \left( \kappa \, v^{2}+\delta^{2} \kappa  - \delta^{2} \lambda \right) + \delta \left[ g' \left( 2\kappa \, v + \delta' \right) - 2 \delta \, \kappa'  + 2\kappa \, g''\,\right]}{2 \delta^{2}\, \kappa \left( \kappa \, g\, v + \delta \, g' \right)}.
    $$
% where the function $g$ satisfies \eqref{420}.
So we arrive at
\medskip
 \begin{proposition}
   A relatively normalized ruled surface  $\varPhi \subset \mathbb{E}^{3}$ is an improper relative sphere \Iff{} the relative normalization is right and one of the following properties holds:
   \begin{description}
     \item[(a)]  $\varPhi$ is conoidal and $f$ and  $g$ are the functions \eqref{405}.
     \item[(b)] $\varPhi$ is non conoidal, the function $g$ fulfils \eqref{420} and $f$ is the function \eqref{415}.
   \end{description}
 \end{proposition}

Let now $\varPhi $ be a proper relative sphere, i.e., by definition \cite{Manhart3}, its relative normals pass through a fixed point.
It is obvious, that this is valid \Iff{} there exists a constant $c\in \rR^{*}$ and a constant vector
$\overline{a} $, such that
    \begin{equation}\label{430}
    \overline{x}=c\,\overline{y}+\overline{a}.
    \end{equation}
Taking into consideration \eqref{125} and \eqref{430}, we observe that $$ B_{i}^{j}=-\frac{\delta_{i}^{j}}{c} \,\,\, \forall \,\,\, i,j=1,2.$$ Consequently, by \eqref{130}, it is
    $$
    H=\frac{-1}{c},\quad K=\frac{1}{c^{2}}.
    $$
Furthermore, taking partial derivatives of \eqref{430} on account of \eqref{10}, \eqref{20}, \eqref{70} and \eqref{76}
% and \eqref{80}
we obtain
\begin{equation}\label{436}
\delta\, \lambda = c\left(y_{1/1}-y_{2}\right), \quad \delta = c \left( \kappa \,y_{2}+y_{3/1} \right),\quad 1 = c \,y_{1/2},\quad 0 = c\,y_{3/2}.
\end{equation}
Because of (\ref{436}d) the support function is again of the form \eqref{320}, i.e., the normalization is right. We distinguish two cases

\textsc{\textbf{Case I}.} $\varPhi$ is conoidal. From \eqref{75} and (\ref{436}c) we  take
    \begin{equation}\label{445}
    g=\frac{c_{3}-\int \!\delta\, \Ud u}{c}, \, c_{3} \in \rR.
    \end{equation}
Then from (\ref{436}a) we have
    \begin{equation*}
    \left( \frac{f}{\delta}\right)^{''}+\frac{f}{\delta}+\frac {\delta\,\lambda}{c}=0,
    \end{equation*}
hence
    \begin{equation}\label{455}
    f=\delta \frac{   \cos u \left( c_{1}+\int\!\delta\,\lambda \sin u \, \Ud u \right ) + \sin u \left ( c_{2}-\int\! \delta \,\lambda \cos u \, \Ud u \right ) }{c}, \, c_{1},c_{2} \in \rR
    \end{equation}
and (\ref{436}b) becomes an identity. For the relative normalization holds \eqref{325}, where the functions $f$ and $g$ are given by \eqref{455} and \eqref{445}, respectively.
   % \begin{equation}\label{460}
       % \begin{split}
      %  \overline{y}=&\frac{v- \cos u  \left(c_{2}-\int\! \delta \,\lambda \cos u \, \Ud u   \right) - \sin u \left( c_{1}+\int\! \delta \,\lambda \sin u \, \Ud u \right)}{c}\overline{e}\\
       % +&\frac{ \cos u \left( c_{1}+\int\!\delta\,\lambda \sin u \, \Ud u \right ) + \sin u \left ( c_{2}-\int\! \delta \,\lambda \cos u \, \Ud u \right ) }{c} \overline{n}-\frac{\int\! \delta \,\Ud u-c_{3}}{c}\overline{z}.
      %  \end{split}
   % \end{equation}
 Conversely, let the relative normalization \eqref{325} be given, where  $f$ and $g$ are the functions \eqref{455} and \eqref{445}, respectively. By using \eqref{10} and \eqref{15} we infer
        \begin{equation}\label{465}
        \overline{s}'=\left(c \, \overline{y}-v\overline{e} \right)'.
        \end{equation}
It follows $$\overline{s}=c \, \overline{y} - v \, \overline{e} + \overline{a},$$ where $\overline{a}$ is a constant vector. Thus, \eqref{430} is valid and therefore $\varPhi$ is a proper relative sphere, whose striction curve $\varGamma$ is parametrized by
    \begin{equation*}
       \begin{split}    \overline{s}=-& \left[ \cos u  \left(c_{2}-\int\! \delta \lambda \cos u \, \Ud u   \right) + \sin u \left( c_{1}+\int\! \delta \lambda \sin u \, \Ud u \right) \right ] \overline{e}- \left(\int\! \delta \Ud u-c_{3}\right)\overline{z}\\
    +& \left [ \cos u \left( c_{1}+\int\!\delta\lambda \sin u  \Ud u \right ) + \sin u \left ( c_{2}-\int\! \delta \lambda \cos u  \Ud u \right ) \right ] \overline{n}+ \overline{a}.
    \end{split}
    \end{equation*}
Finally, from \eqref{210} and \eqref{320} we find
\begin{equation*}
  J=3g \frac{2\delta  g' v- \delta^{2} \lambda g - \delta'f+2\delta  f'}{2\delta^{2} \left( g \, v + f \right)}.
\end{equation*}
%where the functions $f$ and $g$ are given by \eqref{455} and \eqref{445}, respectively.
    %\begin{equation}\label{480}
    %J=\frac{3 \left( \int \!\delta\, \Ud u - c_{3}\right) \Big \{ 2 \cos u \int \!\delta\,\lambda \cos u  \Ud u + \lambda \left( c_{3}-\int \!\delta\, \Ud u \right)+ 2\big[ v- c_{2} \cos u +\sin u \left(c_{1}+\int \delta\lambda \sin u \Ud u\right) \big] \Big \}}{2c \Big\{ v \left(c_{3}- \int \delta \Ud u  \right) +\delta \Big[ \cos u \left( c_{1} + \int \delta\lambda \sin u \Ud u  \right) + \sin u \left( c_{2}-\int \delta\lambda \cos u \Ud u \right) \Big]   \Big\}}\\
    %\end{equation}

 \textsc{\textbf{Case II}.} $\varPhi$ is non conoidal. From \eqref{75} and (\ref{436}c) we  obtain
    \begin{equation}\label{485}
    f=\frac{\delta \left (\delta+ c \, g'\right)}{c \, \kappa}.
    \end{equation}
By using (\ref{436}a) we take
    \begin{equation}\label{486}
    \left[ \frac{\delta \left (\delta+ c \, g'\right)}{c \, \kappa} \right]^{''}+\frac{\delta \left (\delta+ c \, g'\right)}{c \, \kappa}+c\left(\kappa \, g \right)'+\delta \, \lambda=0,
    \end{equation}
and (\ref{436}b) becomes an identity.
The relative normalization results
        \begin{equation}\label{487}
        \overline{y} = \bigg \{ \frac{v}{c}-\kappa \, g -\left[ \frac{\delta \left (\delta+ c \, g'\right)}{c \, \kappa} \right]^{'} \bigg \}\, \overline{e}+\frac{\delta \left (\delta+ cg'\right)}{c\kappa}\,\overline{n} - g \, \overline{z}.
        \end{equation}
Conversely, let the relative normalization \eqref{487} be given, where the function $g$ satisfies \eqref{486}.
If we proceed as in case I, we easily verify that \eqref{465} holds true and therefore $\varPhi$ is a proper relative sphere, whose
%Under the condition \eqref{486} and by using \eqref{15} we obtain
        %\begin{equation}\label{488}         \overline{s}'=\left(c \overline{y}-v\overline{e} \right)'.         \end{equation}
striction curve $\varGamma$ is parametrized by
\begin{equation*}
 \overline{s} = \bigg \{\!\!-c \, \kappa \, g -\left[\frac{\delta \left (\delta + c\,g'\right)}{\kappa}\right]^{'} \bigg \}\,\overline{e}+\frac{\delta \left (\delta+ c \, g'\right)}{\kappa}\,\overline{n}-c\,g\,\overline{z}+\overline{a}.
\end{equation*}
where $\overline{a}$ is a constant vector.
%Thus, \eqref{215} is valid and therefore $\varPhi$ is a proper relative sphere.
Finaly, for the Pick invariant we have
        \begin{equation*}
        J=3g \,\frac{g \left[\kappa \, v^{2} +\delta^{2} (\kappa-\lambda)\right]-\delta' f + 2\delta\left(g'v+f'\, \right)}{2\delta^{2}\left( g \,v + f \right)}.
        \end{equation*}
%where the functions $f$ and $g$ satisfy \eqref{485} and \eqref{486}.
Thus the following has been shown
\medskip
 \begin{proposition}
   A relatively normalized ruled surface  $\varPhi \subset \mathbb{E}^{3}$ is a proper relative sphere \Iff{} the relative normalization is right and one of the following properties holds:
   \begin{description}
     \item[(a)]  $\varPhi$ is conoidal and $f$ and $g$ are the functions \eqref{455} and \eqref{445}, respectively.
     \item[(b)] $\varPhi$ is non conoidal, the function $g$ fulfils \eqref{486} and $f$ is the function \eqref{485}.
   \end{description}
 \end{proposition}

We wish to conclude this section by determining the relative normalizations, which are constantly linked to the moving frame $\mathcal{D}$, i.e., $y_{i|j}=0$ for $i=1,2,3$ and $j=1,2$. An elementary treatment of the last system of equations yields: $\varPhi $ is a conoidal relatively minimal ruled surface and the support function is of the form
%By using \eqref{75} we deduce that the normalization is right, $\varPhi$ is conoidal, $g=c_{1}\in \rR $ and $f=c_{2}\delta$, $c_{2}\in \rR$, $c_{1}^{2}+c_{2}^{2} \neq 0$, i.e., the support function is of the form
\begin{equation*}
  q=\frac{c_{1}\ v + c_{2} \, \delta}{w}, \, c_{1}, c_{2}\in \rR, c_{1}^{2}+c_{2}^{2} \neq 0.
\end{equation*}
Whenever $c_{2}=0$ ($\varPhi, \overline{y}$) is an improper relative sphere ($\overline{y}'=\overline{0}$), otherwise the relative image of $\varPhi$ degenerates into a piece of circle of radius $1 : |c_{2}|$ which is parametrized by
$$\overline{y}=c_{2}\, \overline{n} - c_{1} \, \overline{z}.$$
%\begin{itemize}
 % \item \emph{A relative normalization of a ruled surface $\varPhi \subset \mathbb{E}^{3}$ is constantly linked with the moving frame %$\mathcal{D}$ \Iff{}  $\varPhi$ is conoidal and the support function is of the form \eqref{488}.}
%\end{itemize}

\section{Central normalizations}\label{Sec5}
%According to \cite{Stamatakis6} the relative normal at each point $P$ of $\varPhi$ lies on the corresponding central plane $\{P;\overline{e},\overline{z}\}$ \Iff{} $y_2=0$,or, because of (\ref{75}), \Iff{} the support function of $\overline{y}$ is of the form
 %   \begin{equation}\label{500}
  %  q=\frac{g\,v}{w},
   % \end{equation}
%where $g=g(u)$ is an arbitrary $C^3$-function such that $q>0$. We call the corresponding relative normalization \textit{central}.
%Then from \eqref{70}, \eqref{75} and \eqref{500} it follows
 %   \begin{equation}\label{505}
  %  \overline{y}=-\frac{g'\,v+\delta \,\kappa \, g }{\delta}\overline{e}-g\,\overline{z}.
   % \end{equation}
%The central normalizations are a special case of the right normalizations (\cite{Stamatakis6}).
Let $\overline{y}$ be a central normalization of a given ruled surface $\varPhi$. From \eqref{210} - \eqref{225} and \eqref{301} we obtain
\begin{equation}\label{510}
  H=\frac{g'}{\delta}, \, K=\frac{g'^{2}}{\delta^{2}}, \, S = \frac{-g\left[\kappa v^{2}+\delta^{2}\left(\kappa-\lambda\right)\right]}{2\delta^{2}v}, \, J= \frac{3\left[\kappa g v^{2}+2\delta g' v+ \delta^{2} g \left(\kappa-\lambda\right)\right]}{2\delta^{2}v}.
\end{equation}
It is obvious that $\varPhi$ \emph{is a relative minimal surface} (\emph{or of vanishing relative curvature}) \emph{\Iff{}} $g = c\in\mathbb{R^{*}}$.
% \begin{itemize} is a right conoid} $(\kappa = \lambda = 0)$.
%\begin{itemize}
 % \item $\varPhi$ \emph{is a relative minimal surface (or of vanishing relative curvature) \Iff{}} $g=c\in\mathbb{R^{*}}$.
  %\item \emph{The scalar curvature of the relative metric $G$  vanishes identically \Iff{} $\varPhi$ is a right conoid} $(\kappa = \lambda = 0)$.
  %\item \emph{The Pick invariant vanishes identically \Iff{} $\varPhi$ is a relative minimal right conoid}.
%\end{itemize}
Furthermore, the scalar curvature of the relative metric $G$  vanishes identically \Iff{} $$\kappa v^{2}+\delta^{2}\left(\kappa-\lambda\right) = 0.$$ After successive differentiations of this last equation relative to $v$, we deduce, that $\kappa = \kappa - \lambda = 0$. Thus, \emph{the scalar curvature of the relative metric  vanishes identically \Iff{} $\varPhi$ is a right conoid} $(\kappa = \lambda = 0)$. In the same way one may see that \emph{the Pick invariant vanishes identically \Iff{} $\varPhi$ is a relative minimal right conoid}.

We notice that \emph{all points of $\varPhi $ are relative umbilics} ($H^{2}-K\equiv0 $). Hence, for the relative principal curvatures $k_{1}$ and $k_{2}$ holds $k_{1}=k_{2}=H$. If $g=c\in\mathbb{R^{*}}$, then, because of \eqref{305}, the central normalization of $\varPhi$ reads $\overline{y}=-g\,\overline{d}$, i.e., the central image  degenerates into a curve parallel to the Darboux vector of  $\varPhi$ \cite{Hoschek}. If $g$ is not constant the parametrization of the unique relative focal surface of $\varPhi$ which by definition is given by
\begin{equation*}
\overline{x}^{*}=\overline{s} + v\, \overline{e} + \frac{1}{H} \, \overline{y}
\end{equation*}
becomes
\begin{equation*}
  \overline{x}^{*}=\overline{s} - \frac{\delta \, g}{g'} \,\, \overline{d},
\end{equation*}
i.e., \emph{the focal surface degenerates into a curve $\varGamma ^{* } $ and all relative normals
along each generator form a pencil of straight lines whose centers lie on the curve $\varGamma ^{* }$}.
Let $P(u_{0})$ be a point of the striction curve $\varGamma$ of $\varPhi$ and $R(u_{0})$ the corresponding
point on the focal curve $\varGamma^{*}$. We consider all central normalizations of $\varPhi$. Therefore, the locus of the points $R( u_{0})$ is a straight line parallel to the vector $\overline{d}\left( u_{0}\right) $.
Thus, we obtain a ruled surface $\varPhi ^{* }$, whose generators are parallel to the
vectors $\overline{d}(u)$, a parametrization of which reads
\begin{equation*}
  \overline{x}^{*}=\overline{s} + v^{*} \,\, \overline{d}.
\end{equation*}
A parametrization of its striction curve is
\begin{equation*}
  \overline{s}^{*}=\overline{s} - \frac{\delta \, \lambda}{\kappa'} \, \overline{d}.
\end{equation*}
One can easily verify that $\varPhi^{*}$ is developable.

By using \eqref{305} and \eqref{400} we may infer:
\begin{itemize}
  \item \emph{A centrally normalized ruled surface $\varPhi$ is an improper relative sphere \Iff{}
  %$g=c_{1}\in\mathbb{R^{*}}$ and $\kappa=c_{2}\in\mathbb{R}$, i.e.
   $\varPhi$ is a relative minimal surface of constant conical curvature}. Then, the central image of $\varPhi$ degenerates into a curve parallel to the Darboux vector.
   \item \emph{A centrally normalized ruled surface $\varPhi$ is a proper relative sphere \Iff{}}
\begin{equation*}
  g = \frac{1}{c}\left(c_{1}-\int \!\! \delta \, \Ud u \right), c, c_{1}\in\mathbb{R} \quad \text{and} \quad \delta(\kappa-\lambda)+\kappa'\left(\int \!\! \delta \, \Ud u -c \, c_{1} \right)=0.
\end{equation*}
 \end{itemize}
We focus now on the field $\overline{T}(u,v)$ of the Tchebychev vectors of $(\varPhi, \overline{y}$). By using \eqref{120} we find
\begin{equation}\label{705}
  \overline{T}=\frac{2 \kappa g v^{2}+\left(\delta' g+2\delta g '\right)v+2\kappa\delta^{2}g}{2\delta^{2}}\overline{e}+\frac{g}{\delta}\left(v \overline{n}+\delta\overline{z}\right).
\end{equation}
By taking (\ref{240}) and (\ref{510}c) into consideration it turns out that the divergence $\divz ^{G}\overline{T}$ of $\overline{T}$ with respect to the relative metric $G$ of $\varPhi$ is given by
\begin{equation*}
\divz ^{G}\overline{T}=\frac{g}{\delta^{2}v}\left[\kappa v^{2}-\delta^{2}(\kappa-\lambda)\right].
\end{equation*}
Consequently we have $$\divz ^{G}\overline{T}  = 0$$ \Iff{} $$\kappa v^{2} - \delta^{2}(\kappa - \lambda) = 0.$$ After successive differentiations of this last equation relative to $v$, we deduce, that $\kappa = \kappa - \lambda = 0$.
% \begin{itemize}
  %\item
So we have:  \emph{The vector field} $\overline{T}$ \emph{is incopressible \Iff{}} $\varPhi$ \emph{is a right conoid.}
%\end{itemize}
%Our aim is to classify the central normalized ruled surfaces $\varPhi$, whose Tchebychev vector field is tangent or orthogonal to some geometrically distinguished families of curves of $\varPhi$.

The vectors $\overline{T}$ are orthogonal to the generators of $\varPhi$ \Iff{} $$\langle \overline{e},\overline{T}\rangle=0.$$ On account of \eqref{705} we have
\begin{equation*}
  2 \kappa g v^{2}+\left(\delta'g+2\delta g'\right)v+2\delta^{2}\kappa g=0.
\end{equation*}
%On differetiating three times relative to $v$ we infer the system
%\begin{equation*}
 % 2 \kappa g=\delta'g+2\delta g'=2 \delta^{2} \kappa g= 0,
%\end{equation*}
Treating analogously this equation we conclude: $\overline{T}$ \emph{is orthogonal to the generators of} $\varPhi$ \emph{\Iff{}} $\varPhi$ \emph{is conoidal and} $g = c |\delta|^{-1/2},\, c\in\mathbb{R^{*}}$.
%which implies that $\kappa=0$ and $g= c|\delta|^{-1/2}, c\in\mathbb{R^{*}}$. The inverse also holds. %So we obtain:
%\begin{itemize}
 % \item \emph{$ \overline{T}$ is orthogonal to the generators of $\varPhi$ \Iff{} $\varPhi$ is conoidal and $g= c|\delta|^{-1/2}, c\in\mathbb{R^{*}}$}.
%\end{itemize}

The tangent vector to a curve $ \varLambda : v=v(u)$ of $\varPhi$ is
\begin{equation}\label{730}
  \overline{x}'=\left(\delta \, \lambda +v '\right)\,\overline{e} + v \, \overline{n} + \delta \,\overline{z}.
\end{equation}
We consider now the following families of curves on $\varPhi$: a) the curved asymptotic lines, b) the curves of constant striction distance ($u$-curves) and c)  the $\widetilde{K}$-curves, i.e., the curves along which the Gaussian curvature is constant \cite{Sachs}.
%\begin{description}
 % \item[a.] the curved asymptotic lines,
  %\item[b.] the curves of constant striction distance ($u$-curves) and
  %\item[c.] the $\widetilde{K}$-curves, i.e., the curves along which the Gaussian curvature is constant \cite{Sachs}.
%\end{description}
The corresponding differential equations of these families of curves are
\begin{align}
  &\kappa v^{2}+\delta'v+\delta^{2}\left(\kappa-\lambda\right)-2\delta v'=0, \label{715} \\
  &v'=0, \\
  & 2 \delta v v'+ \delta'\left(\delta^{2}-v^{2}\right)=0.
\end{align}
From \eqref{705} and \eqref{730} it follows: $ \overline{x}'\,$ and $\overline{T}$ are parallel or orthogonal \Iff{}
\begin{equation}\label{735}
2 \kappa g v^{2} + \left(\delta' g+ 2 \delta g'\right)v+ 2 \delta^{2} \kappa g-2 \delta g \left(\delta\lambda+v'\right)=0
\end{equation}
or
\begin{equation}
  \left(\delta\lambda+v'\right)\left[2 \kappa g v^{2} + \left(\delta' g+ 2 \delta g'\right)v+ 2 \delta^{2} \kappa g \right] + 2\delta g w^{2}=0,
\end{equation}
respectively. From \eqref{715} and \eqref{735} we infer, that $\overline{T}$ \emph{is tangential to the curved asymptotic lines \Iff{}}
\begin{equation*}
  \kappa g v^{2}+ 2 \delta g' v + \delta^{2} g \left(\kappa-\lambda\right)=0,
\end{equation*}
that is \Iff{} $\kappa=\lambda=0$ and $g=const.$ and therefore \emph{\Iff{}} $\varPhi$ \emph{is a relative minimal right conoid}.
%\begin{itemize}
 % \item\emph{$ \overline{T}$ is tangential to the curved asymptotic lines of $\varPhi$ \Iff{} $\varPhi$ is a relative minimal right conoid.}
%\end{itemize}

Analogous reasoning in the case of the families of curves b) and c) leads to the following results: %From \eqref{720} and \eqref{735}, resp. \eqref{740}, we have: $\overline{T}$ is tangential or orthogonal to
% the $u$-curves \Iff{}
% \begin{equation*}
%   2 \kappa g v^{2}+\left(\delta' g + 2 \delta g'\right)v+ 2 \delta^{2}g(\kappa-\lambda)=0
% \end{equation*}
% or
% \begin{equation*}
%   2g (\kappa\lambda+1)v^{2}+\lambda (\delta' g + 2 \delta g')v + 2\delta^{2}g (\kappa\lambda+1)=0,
% \end{equation*}
% respectively. Treating analogously these polynomials we find:
 \begin{itemize}
   \item \emph{$ \overline{T}$ is tangential to the $u$-curves of $\varPhi$ \Iff{} $\varPhi$ is right conoid and $g= c|\delta|^{-1/2}$, $c\in\mathbb{R^{*}}$.}
   \item\emph{$ \overline{T}$ is orthogonal to the $u$-curves of $\varPhi$ \Iff{} the striction curve of $\varPhi$ is an Euclidean line of curvature and $g= c|\delta|^{-1/2}, c\in\mathbb{R^{*}}$.}
 \end{itemize}
%From \eqref{725} and \eqref{735}, resp. \eqref{740}, we infer: $\overline{T}$ is tangential or orthogonal to
% the $\widetilde{K}$-curves \Iff{}
% \begin{equation*}
%   2 \kappa g v^{3}+ 2 \delta g' v^{2}+2 \delta^{2}g (\kappa-\lambda) v+\delta^{2}\delta' g=0
% \end{equation*}
% or
% \begin{equation*}
%   \begin{split}
%    & 2 \delta' \kappa g v^{4}+ \left[ 4\delta^{2}g(\kappa\lambda+1)+\delta'(\delta' g+ 2 \delta g')\right]v^{3}+2 \delta^{2}\lambda \left(\delta' g +2\delta g'\right)v^{2} \\
 %   &+\delta^{2}\left[4\delta^{2}g(\kappa\lambda+1)-\delta'(\delta'g+2\delta g')\right]v-2 \delta^{4}\delta'\kappa g=0,
%   \end{split}
% \end{equation*}
%respectively. Treating these polynomials in the same way we arrive at:
\begin{itemize}
  \item\emph{$ \overline{T}$ is tangential to the $\widetilde{K}$-curves of $\varPhi$ \Iff{} $\varPhi$ is a relative minimal right helicoid surface} $(\kappa = \lambda = 0, \delta = const.)$.
  \item\emph{$ \overline{T}$ is orthogonal to the $\widetilde{K}$-curves of $\varPhi$ \Iff{} $\varPhi$ is a relative minimal Edlinger surface} ($\kappa \, \lambda + 1 = 0, \delta = const.$, see\cite{Hoschek}, \cite{Sachs}).
\end{itemize}
Next we assume that the vector field $\overline{T}$ is tangential to one family of Euclidean lines of curvature. Their differential equation which initially reads
\begin{equation*}
  g_{12}h_{11}-g_{11}h_{12}+\left(g_{22}h_{11}-g_{11}h_{22}\right)v'\,+(g_{22}h_{12}-g_{12}h_{22})v'^{2}=0,
\end{equation*}
becomes, on account of \eqref{30} and \eqref{35},
\begin{equation*}
  \delta \left[w^{2}\left(\kappa \lambda +1\right)+ \delta' \lambda v\right]+\left(\kappa w^{2}+\delta'v-\delta^{2}\lambda\right) v'-\delta v'^{2}=0,
\end{equation*}
from which, by virtue of \eqref{735}, we infer %, that $\overline{T}$ is tangential to one family of the lines of curvature of $\varPhi$ \Iff{}
\begin{equation*}
  \begin{split}
      & 2 \kappa g \left(\delta' g-2 \delta g'\right)v^{3} +\left[4 \delta^{2} g^{2} \left(\kappa \lambda+1\right)+\delta'^{2}g^{2}-4\delta^{2} g'^{2}\right]v^{2} \\
       & + 2 \delta^{2}g \left[\delta'g(\kappa+\lambda)-2\delta g'(\kappa-\lambda)\right]v+4\delta^{4}g^{2}\left(\kappa \lambda+1\right)=0.
  \end{split}
\end{equation*}
This last equation holds true \Iff{} $\kappa \, \lambda + 1 = 0, \delta = c_1 \in \rR^* $ and $ g = c_2 \in \rR^*$.
Hence  \emph{$\varPhi$ is a relative minimal Edlinger surface}. Because of \eqref{15} and
\begin{equation*}
  \overline{T}_{v=0}=g\left(\kappa\overline{e}+\overline{z}\right)
\end{equation*}
we obtain $\langle\overline{T}_{v=0}, \overline{s}'\,\rangle=0 $. So, we have:
%\begin{itemize}
 % \item $ \overline{T}$ \emph{ is tangential to the family of the lines of curvature of $\varPhi$ to which the striction curve does not belong \Iff{} $\varPhi$ is a relative minimal Edlinger surface}.
%\end{itemize}
\medskip
\begin{proposition}
$ \overline{T}$ is tangential to the one family of the Euclidean lines of curvature of $\varPhi$ \Iff{} $\varPhi$ is a relative minimal Edlinger surface. Moreover the family of the Euclidean lines of curvature under consideration consists of the lines of curvature which are orthogonal to the striction curve.
\end{proposition}
In the rest of this section se assume that $\varPhi$ is a non relatively minimal ruled surface. The central image $\varPsi_{1}$ of $\varPhi$ is also a ruled surface, whose generators are parallel to those of $\varPhi$. Then, from \eqref{305} by direct computation, we find that the parametrization of its striction curve is
\begin{equation*}
\text{$\varGamma_{1}$:} \quad  \overline{s}_{1} = -g \, \overline{d}.
\end{equation*}
We set $\overline{y}=\overline{y}_{1}$ and we rewrite the parametrization of $\varPsi_{1}$ as
\begin{equation*}
  \text{$\varPsi_{1}$:} \quad \overline{y}_{1} = \overline{s}_{1} + v_{1} \, \overline{e}, \quad v_{1}:=- \frac{g' v}{\delta}.
\end{equation*}
It is obvious that $\varPsi_{1}$ is parametrized like in (\ref{1}) and (\ref{5}). We use $\mathcal{D}$ as moving frame of $\varPsi_{1}$. By direct computation we find the fundamental invariants of $\varPsi_{1}$:
\begin{equation}\label{610}
  \kappa_{1}=\kappa, \quad \delta_{1}=-g', \quad \lambda_{1}=\frac{\left(\kappa \, g\right)^{'}}{g'}.
\end{equation}
Thus the Darboux vectors of $\varPhi$  and $\varPsi_{1}$ are parallel. Furthermore it is %$w_{1}^{2}=v_{1}^{2}+\delta_{1}^{2}=\frac{g'^{2}w^{2}}{\delta^{2}}$ or equivalently
$w_{1}=|H| w$. Then, by (\ref{40}), the Gaussian curvature $\widetilde{K}_{1}$ of $\varPsi_{1}$ is seen to be
\begin{equation*}
  \widetilde{K}_{1}=\frac{\widetilde{K}}{K}.
\end{equation*}
From the above we list the following results, which can be checked easily:
\begin{itemize}
  \item\emph{$\varPhi$  and its central image $\varPsi_{1}$ are congruent $(\delta=\delta_{1}, \kappa=\kappa_{1}, \lambda=\lambda_{1})$ \Iff{} $\varPhi$ is a proper relative sphere.}
  \item \emph{$\varPsi_{1}$ is orthoid $(\lambda_{1}=0)$ \Iff{} $\kappa g = c \in\mathbb{R}$.}
  \item \emph{The striction curve $\varGamma_{1}$ of $\varPsi_{1}$ is an asymptotic line of it $(\kappa_{1}=\lambda_{1})$ \Iff{} $\varPhi$ is of constant conical curvature}.
  \item \emph{The striction curve $\varGamma_{1}$ of $\varPsi_{1}$ is an Euclidean line of curvature of it $( 1 + \kappa_{1} \, \lambda_{1}=0)$ \Iff{}}
  \begin{equation*}
    g=\frac{c}{\sqrt{1+\kappa^{2}}},\, c\in\mathbb{R^{*}}.
  \end{equation*}
  \item \emph{$\varPsi_{1}$ is an Edlinger surface \Iff{}}
  \begin{equation*}
    g=c_{1}+c_{2}u, \, c_{1}\in\mathbb{R}, \, c_{2}\in\mathbb{R^{*}} \quad \text{and} \quad 1+\kappa^{2}=\frac{c^{2}}{\left(c_{1}+c_{2}u\right)^{2}}, \, c\in\mathbb{R^{*}}.
  \end{equation*}
\end{itemize}

We wish to conclude this paper by answering the following question: Is there a ruled surface $\varPsi^{*}$ whose a central normalization is the given ruled surface $\varPhi$ ? We suppose that such a ruled surface exists and let it be parametrized like in (\ref{1}) and (\ref{5}). We consider a central normalization of $\varPsi^{*}$ via a support function of the form \eqref{301}. Denoting by  $\delta^{*}$, $ \kappa^{*}$ and $\lambda^{*}$ its fundamental invariants we have:
\begin{equation*}
  \kappa = \kappa^{*}, \quad \delta = -g^{*}\,', \quad \lambda = \frac{\left( \kappa \, g^{*} \right)^{'}}{g^{*}\,'},
\end{equation*}
cf. \eqref{610}, thus
\begin{equation}\label{611}
g^{*} = c_{1}- \int \!\! \delta \, \Ud u \,\, \left(c_{1}\in\mathbb{R}\right), \quad \delta \left( \kappa - \lambda \right) - \kappa' \left( c_{1} - \int \!\! \delta \, \Ud u \right) = 0.
\end{equation}

 Conversely, we assume that a constant $c_{1}\in\mathbb{R}$ exists, such that the fundamental invariants of $\varPhi$ fulfil (\ref{611}b). We consider an arbitrary skew ruled surface $\varPsi^{*}$, whose generators are parallel to those of $\varPhi$ and we normalize it centrally via a support function of the form \eqref{301}, where $g^{*}$ is a function of the form (\ref{611}a). By taking \eqref{610} and (\ref{611}b) into account, we deduce that the fundamental invariants of the arising central image of $\varPsi^{*}$  are $\delta$, $\kappa$ and $\lambda$. Hence, the central image of $\varPsi^{*}$ and $\varPhi$ are congruent. Summing up we have:
\medskip
\begin{proposition}
A necessary and sufficient condition for the existence of a ruled surface $\varPsi^{*}$, whose central normalization is the given ruled surface $\varPhi$ is that there is a constant $c_{1}\in\mathbb{R}$, such that the fundamental invariants of $\varPhi$ fulfil $(\ref{611}b)$.
\end{proposition}
An example of a ruled surface, whose fundamental invariants fulfil $(\ref{611}b)$ is the right helicoid.

\end{document}